\documentclass[pdftex,12pt,a4paper,reqno]{amsart}
\usepackage[T1]{fontenc}
\usepackage[latin1]{inputenc}
\usepackage[english]{babel}
\usepackage{lmodern}
\usepackage{amsthm}
\usepackage{amsmath}
\usepackage{amsfonts}
\usepackage{amssymb}
\usepackage[all]{xy}
\usepackage{calrsfs}
\usepackage{slashed}
\usepackage{babel}
\usepackage{array}
\usepackage{tensor}
\theoremstyle{plain}
\newtheorem{theo}{Theorem}[section]
\newtheorem{prop}[theo]{Proposition}
\newtheorem{lemm}[theo]{Lemma}

\newtheorem{cor}[theo]{Corollary}

\theoremstyle{definition}
\newtheorem{defi}[theo]{Definition}

\theoremstyle{remark}

\newtheorem{ex}[theo]{Example}

\newtheorem{rk}[theo]{Remark}

\newcommand{\C}{\mathbb{C}}

\newcommand{\Heis}{\mathbb{H}}
\newcommand{\N}{\mathbb{N}}

\newcommand{\R}{\mathbb{R}}

\newcommand{\Sph}{\mathbb{S}}

\newcommand{\re}{\Re\mathrm{e}}

\usepackage[bookmarks,bookmarksopen,bookmarksdepth=2]{hyperref}
\usepackage[usestackEOL]{stackengine}
\usepackage{etoolbox}
\patchcmd{\thmhead}{(#3)}{#3}{}{}
\usepackage[left=3cm,right=3cm,top=3cm,bottom=3cm]{geometry}

\everymath{\displaystyle\everymath{}}

\begin{document}
\title{Contact structures, CR Yamabe invariant, and connected sum}
\author{Gautier Dietrich}
\thanks{The author was supported in part by the grant ANR-17-CE40-0034 of the French National Research Agency ANR (project CCEM)}
\address{Institut Montpelliérain Alexander Grothendieck\\ Université de Montpellier\\ CNRS\\ Case courrier 051\\ Place Eugène Bataillon\\ 34090 Montpellier\\ France} 
\address{Université Paul-Valéry Montpellier 3%\\ UFR 6\\ Route de Mende \\34199 Montpellier cedex 5\\ France
} 
\email{gautier.dietrich@ac-montpellier.fr}
\date{}

\maketitle

\begin{abstract}
We propose a global invariant $\sigma_c$ for contact manifolds which admit a strictly pseudoconvex CR structure, analogous to the Yamabe invariant $\sigma$. We prove that this invariant is non-decreasing under handle attaching and under connected sum. We then give a lower bound on $\sigma_c$ in a particular case.
\end{abstract}

\setcounter{tocdepth}{1}
\tableofcontents

\section{Introduction}

%[Contact invariants are scarce, and identifying them is an important issue. Note that from Darboux's theorem, contact invariants are necessarily global.]
The classical \emph{Yamabe problem} is that of the existence, on a compact Riemannian manifold, of a metric conformal with a given metric and with constant scalar curvature \cite{Sch84,LP87}. The \emph{Yamabe invariant} $\sigma$ has been introduced by R. Schoen and O. Kobayashi in the wake of the resolution of the Yamabe problem \cite{Kob87,Sch89}. It is built the following way: a sufficient condition for a metric $g$ to have constant scalar curvature is to minimize, among metrics of same volume in the same conformal class, the integral scalar curvature $S(g)$. This minimum is moreover always smaller than $S(g_{\Sph^n})$, where $g_{\Sph^n}$ is the standard metric on the sphere. The Yamabe invariant is then defined, for a compact differentiable manifold $M$, as the "max-min" $$\sigma(M):=\sup_{[g]}\inf_{\hat{g}\in [g]_1} \int_M\mathrm{Scal}(\hat{g})\ d\text{vol}_{\hat{g}},$$ where $\mathrm{Scal}$ denotes the Riemannian scalar curvature, the supremum runs over all conformal classes of metrics $[g]$ on $M$ and the infimum runs over all metrics of volume $1$ in $[g]$.

This global differential invariant is produced by looking at the given differentiable manifold $M$ through the prism of the conformal structures $[g]$ with which it can be equipped. Similarly, let us consider a compact contact manifold $(M,H)$ which admits a strictly pseudoconvex CR structure, which we will call an \emph{SPC} manifold. To each CR structure corresponds a conformal class of positive contact forms $\theta$, which leads to the \emph{CR Yamabe problem}. %Let us consider the 
This problem has been given a positive answer by D. Jerison and J. Lee, N. Gamara and R. Yacoub \cite{JL87,JL89,Gam01,GY01}. One can thus define a \emph{contact Yamabe invariant} $\sigma_c$ as% for contact manifolds $(M,H)$ that can be equipped with an SPCR structure, or \textit{SPC} manifolds \cite{Wu09},
 $$\sigma_c(M,H):=\sup_{\mathcal{J}} \inf_{\hat{\theta}\in [\theta]_1}\int_M\mathrm{Scal}_W(J,\hat{\theta})\ \hat{\theta}\wedge (d\hat{\theta})^n,$$ where $\mathrm{Scal}_W$ denotes the \emph{Webster} scalar curvature, the supremum runs over the set $\mathcal{J}$ of all complex structures $J$ such that $(M,H,J)$ is strictly pseudoconvex, and the infimum runs over all compatible pseudohermitian forms of volume 1%, and $\int_{(M,\hat{\theta})}\mathrm{Scal}_W$ is the integral of the Webster scalar curvature over $M$ with respect to $\hat{\theta}$
. This invariant has been introduced by C.-T. Wu \cite{Wu09}, but up to our knowledge has not been studied since. This invariant is a genuine contact invariant in dimension $3$, in the sense that all orientable $3$-dimensional contact manifolds are SPC. %This is no longer the case in higher dimension.

A construction by W. Wang, recently implemented by J.-H. Cheng and H.-L. Chiu, shows that the positivity of $Y_{CR}$ is preserved under handle attaching on a strictly pseudoconvex CR spherical manifold \cite{Wan03,CC19}:

\begin{theo}[\cite{Wan03,CC19}]
Let $(M,H,J)$ be a compact spherical strictly pseudoconvex CR manifold with positive $Y_{CR}$. Let $(\tilde{M},\tilde{H},\tilde{J})$ be obtained from $(M,H,J)$ by CR handle attaching. Then $(\tilde{M},\tilde{H},\tilde{J})$ is spherical and $Y_{CR}(\tilde{M},\tilde{H},\tilde{J})>0$.
\end{theo}

From a continuity argument detailed in Section \ref{lsph}, we generalize this result:

\begin{theo}
\label{CRkobh}
Let $(M,H)$ be a compact SPC manifold. Let $(\tilde{M},\tilde{H})$ be a manifold obtained from $(M,H)$ by SPC handle attaching, then $$\sigma_c(\tilde{M},\tilde{H})\geq \sigma_c(M,H).$$
\end{theo}

\noindent Moreover, we prove the contact analogue of a theorem due to Kobayashi \cite{Kob87}:

\begin{theo}
\label{CRkob}
Let $(M_1,H_1)$ and $(M_2,H_2)$ be two compact SPC manifolds of dimension $2n+1$. Let $(M_1,H_1)\# (M_2,H_2)$ be their SPC connected sum, then $$\sigma_c\left((M_1,H_1)\# (M_2,H_2)\right)\geq \left\lbrace
\begin{array}{ll}
-\left(|\sigma_c(M_1,H_1)|^{n+1}+|\sigma_c(M_2,H_2)|^{n+1}\right)^\frac{1}{n+1} & \mathrm{if} \ \sigma_c(M_1,H_1)\leq 0\\
& \hspace{10px} \mathrm{and} \ \sigma_c(M_2,H_2)\leq 0,\\
\min\left(\sigma_c(M_1,H_1),\sigma_c(M_2,H_2)\right) & \mathrm{otherwise}.
\end{array}\right.$$
\end{theo}

%In this paper, we generalize this construction to prove a CR analogue of Theorem \ref{adh}:
%\begin{theo}
%\label{CRineq}
%Let $(M,H)$ be a compact SPC manifold. When $n=k=1$, there is a positive constant $\Lambda^{CR}_{n,k}$ such that for all $(N,\tilde{H})$ obtained from $(M,H)$ by a k-dimensional Eliashberg surgery, $$\sigma_c(N,\tilde{H})\geq\min\left(\sigma_c(M,H),\Lambda^{CR}_{n,k}\right).$$
%\end{theo}

We also prove a weakened contact version of a theorem due to C. LeBrun and J. Petean, who, using the generalized Gauss-Bonnet theorem, have computed the Yamabe invariant for complex surfaces %(of real dimension $4$) 
of general type \cite{LeB96,Pet98}:

%\noindent Our result gives the value of the $\sigma_c$ invariant in the following case:
%
%\begin{theo}
%%Let $(X,H)$ be a circle bundle over a Riemann surface $\Sigma$ of positive genus. If, for all $J\in\mathcal{J}$, the unique normal contact form on $(X,H,J)$ is also pseudo-Einstein %if $\mathcal{J}\ni J\mapsto \mu(M,H,J)$ realizes its supremum at $J_0$, 
%%and %there exists a Yamabe contact form
%%Yamabe, then $$\sigma_c(X,H)=-2\pi\sqrt{-\chi(\Sigma)},$$ where $\chi(\Sigma)$ is the Euler characteristic of $\Sigma$.
%Let $(M,H)$ be an SPC compact $3$-manifold admitting an Einstein, Yamabe pseudohermitian structure $(J,\theta)$ with nonpositive Burns-Epstein invariant $\mu(M,H,J)$. Then $$\sigma_c(M,H)\geq \varepsilon 4\pi%\underset{J\in\mathcal{J}}{\sup}
%\sqrt{-\mu(M,H,J)},$$ where $\varepsilon=\mathrm{sign}\left(Y_{CR}(M,H,J)\right).$
%\end{theo}

\begin{theo}\label{CRleb}
Let $(M,H)$ be a circle bundle over a Riemann surface $\Sigma$ of positive genus admitting an Einstein pseudohermitian structure. Then $$\sigma_c(M,H)\geq -2\pi\sqrt{-\chi(\Sigma)}.$$
\end{theo}

%Section \ref{CRdef} is devoted to notions of CR geometry. Section \ref{CRYam} is an overview of the CR Yamabe problem.% Section \ref{CRsur} details a contact surgery construction preserving the strict pseudoconvexity, after \cite{Eli90}. 
Section \ref{CRadh} contains the proof of Theorems \ref{CRkobh} and \ref{CRkob}, and Section \ref{comp} contains the proof of Theorem \ref{CRleb}.

\textbf{Acknowledgements.} I am deeply grateful to my supervisor, Marc Herzlich, for introducing me to these questions, for his numerous advices and his precious help. I also thank the referees of my PhD thesis, Jih-Hsin Cheng and Colin Guillarmou, for their careful reading and for pointing out some mistakes. In particular, the proof of Lemma \ref{locsphe2} has greatly benefited from the help of Sylvain Brochard and from the work of Jih-Hsin Cheng, Hung-Lin Chiu, and Pak Tung Ho \cite{CCH19}.

\section{CR geometry}
\label{CRdef}

\subsection{Generalities}

Let $n\in\N^*$ and $M$ be a smooth differentiable manifold of real dimension $2n+1$. We assume that $M$ is orientable. A \emph{CR structure} is given on $M$ by a complex subbundle $T^{1,0}M$ of $TM\otimes \C$ of complex dimension $n$ verifying $$T^{1,0}M\cap T^{0,1}M=\lbrace 0\rbrace$$ where $T^{0,1}M=\overline{T^{1,0}M}$, and which is stable under the Lie bracket.

Equivalently, let $H$ be a \emph{Levi distribution}, \emph{i.e.} an orientable hyperplane distribution in $TM$. Let $J$ be a \emph{complex structure} on $H$, \emph{i.e.} $J$ is an endomorphism of $H$ which satisfies $J^2=-\mathrm{id}_H$ and is \emph{integrable}: $\forall X, Y \in \Gamma(H),$ $$[JX,Y]+[X,JY]\in\Gamma(H)\quad \text{and}\quad [JX,JY]-[X,Y]=J\left([JX,Y]+[X,JY]\right),$$ where $[\cdot,\cdot]$ denotes the Lie bracket. The existence of $J$ also requires that $H$ is orientable. A \emph{CR manifold} is the triplet $(M,H,J)$.

%The correspondence between the two approaches is the following:
%
%\begin{itemize}
%
%\item[$\bullet$] From the first approach to the second, $H:=\re\left(T^{1,0}M\oplus T^{0,1}M\right)$, and $\forall X\in T^{1,0}M,\ J(X+\overline{X}):=i(X-\overline{X})$.
%
%\item[$\bullet$] From the second approach to the first, $T^{1,0}M$ is the $i$-eigenspace of the extension of $J$ to $H\otimes\C$.
%
%\end{itemize}

%An \emph{almost CR structure} on $M$ is a complex subbundle $T_{1,0}M$ of $TM\otimes \C$ of complex dimension $n$ verifying $$T_{1,0}M\cap\overline{T_{1,0}M}=\lbrace 0\rbrace.$$ This structure is called \emph{integrable}, or a \emph{CR structure}, if $\Gamma(T_{1,0}M)$ is stable under the Lie bracket.
%
%Let us denote $T_{0,1}M:=\overline{T_{1,0}M}$, $\C H:=T_{1,0}M\oplus T_{0,1}M$ and $J\in \mathrm{End}(\C H)$ such that $J\mid_{T_{1,0}M}=+i$, $J\mid_{T_{0,1}M}=-i$. The \emph{Levi distribution} of a CR structure $T_{1,0}M$ is the hyperplane distribution $H$ in $TM$ given by $H:=\mathrm{Re}(\C H)$. The \emph{complex structure} $J$ is well-defined on $H$: for all $X\in H$ there is a unique $Z\in T_{1,0}M$ such that $X=Z+\bar{Z}$, and then $JX=i(Z-\bar{Z})$. A \emph{CR manifold} is the triplet $(M,H,J)$.

Let $E:=\lbrace \omega\in \Gamma(T^* M) \ | \ \ker\omega \supseteq H\rbrace\simeq TM/H$. It is a real line subbundle of $T^* M$, hence trivial since $M$ is orientable. A \emph{pseudohermitian structure} on $M$ is a never-vanishing section $\theta$ of $E$ compatible with $J$, \emph{i.e.} such that $$d\theta(J\cdot,J\cdot)=d\theta(\cdot,\cdot)\quad \text{on}\ TM\otimes\C.$$ The associated \emph{Levi form} $\gamma$ is the Hermitian form on $H$ given by $\gamma:=d\theta(\cdot,J\cdot)$.

\begin{defi}
A pseudohermitian structure $\theta$ is said to be \emph{strictly pseudoconvex} when its Levi form is definite positive and when the orientation of the associated volume form $\theta\wedge d\theta^n$ coincides with the orientation of $M$.
\end{defi}

In that case, $\theta$ is a contact form, and $(M,H)$ is a contact manifold.  A contact form on $(M,H,J)$ which is a strictly pseudoconvex pseudohermitian structure will be called \emph{positive}. A CR manifold admitting an positive contact form is called \emph{SPCR}, and a contact manifold admitting an SPCR structure is called \emph{SPC}. We will always assume that $H$ is a contact distribution.% A contact manifold is an SPC manifold if and only if its \emph{Tanno tensor} vanishes \cite{Tan89}. 

\begin{defi}
Given an SPC manifold $(M,H)$, we define $$\mathcal{J}=\left\lbrace \text{Compatible complex structures } J \text{ on }H\ |\ (M,H,J)\text{ is an SPCR manifold}\right\rbrace.$$
\end{defi}

In dimension $2n+1=3$, $T^{1,0}M$ is of complex rank $1$, the integrability of $J$ is thus automatic. The set $\mathcal{J}$ is defined by purely algebraic conditions, and it is moreover contractible. Indeed, considering $J_0,\ J_1$ in $\mathcal{J}$, let, for $i$ in $\lbrace 0,1\rbrace$, $\gamma_i:=d\theta(\cdot,J_i\cdot)$. For $t$ in $[0,1]$, the metric $\gamma_t:=(1-t)\gamma_0+t\gamma_1$ gives a complex structure $\tilde{J}_t$ compatible with $\theta$, with $\tilde{J}_0=J_0$ and $\tilde{J}_1=J_1$. The set $\mathcal{J}$ is therefore always non-empty if $M$ is orientable.

%\begin{defi}
%A CR manifold $(M,H,J)$ is called \emph{embeddable} if there exists a submanifold $(M',H',J')$ of $\C^N$ for some $N$ 	and a CR isomorphism $f:M\rightarrow M'$, \emph{i.e.} such that $f_*(T^{1,0}M)= T^{1,0}M'$.
%\end{defi}
%\noindent All compact CR manifolds of dimension $\geq 5$ are embeddable, but there exist CR $3$-manifolds, even close to $\Sph^3$, that are non-embeddable \cite{Ros65,Bou75}.
%
%As already mentioned, an SPCR structure induces a conformal structure on the contact distribution $H$. Indeed, let $(M,H,J)$ be an SPCR manifold and $\theta$ be a positive contact form on $M$. Note that $E$ is the set of sections of a line bundle. Strict pseudoconvexity is then preserved under conformal changes of $\theta$. Indeed, $$\forall f\in C^\infty(M,\R_+^*),\quad\gamma_{f\theta}=f \gamma_\theta.$$

The \emph{Reeb field} of a contact form $\theta$ is the unique vector field $R\in TM$ verifying $\theta(R)=1$ and $\iota_R d\theta=0$. We get a pseudohermitian decomposition of the tangent space $$TM=\R R\oplus H,$$ and a pseudohermitian projection $\pi_b:TM\rightarrow H$. Note that this projection depends on $\theta$. An \emph{admissible coframe} is a set of $(1,0)$-forms $(\theta_1,\ldots,\theta_n)$ whose restriction to $T^{1,0}M$ forms a basis for $\left(T^{1,0}M\right)^*$ and such that, for all $\alpha$ in $\lbrace 1,\ldots,n\rbrace$, $\theta^\alpha(R)=0$. Then $d\theta=i h_{\alpha\overline{\beta}}\theta^\alpha\wedge\theta^{\overline{\beta}}$, where $\theta^{\overline{\beta}}:=\overline{\theta^\beta}$ and $(h_{\alpha\overline{\beta}})$ is a positive definite Hermitian matrix. If $(T_1,\ldots,T_n)$ is the dual frame to $(\theta^1,\ldots,\theta^n)$ on $T^{1,0}M$, then $\gamma(T_\alpha,T_{\overline{\beta}})=h_{\alpha\overline{\beta}}.$

%\begin{defi}
%The Riemannian metric $g_{J,\theta}$ on $M$ given by $$g_{J,\theta}:=\theta^2+\gamma,$$ is called the \emph{Webster metric} of $(M,H,J,\theta)$.
%\end{defi}

\begin{prop}[\cite{Tan75,Web78}]\label{tw}
Let $(M,H,J,\theta)$ be a strictly pseudoconvex pseudohermitian manifold. There is a unique linear connection $\nabla^\theta$ on $M$, called the \emph{Tanaka-Webster connection}, which parallelizes the Levi distribution $H$, the Reeb field $R$, the complex structure $J$, and the Webster metric $g_{J,\theta}$, and whose torsion $T^\theta$ verifies $$\forall X,Y\in H,\quad T^\theta(X,Y)=d\theta(X,Y)R \quad \mathrm{and} \quad T^\theta(R,JX)+JT^\theta(R,X)=0.$$
\end{prop}
%
%\begin{prop}
%$\nabla^\theta R=0$.
%\end{prop}
\noindent In other words, if $\theta^1,\ldots,\theta^n$ is an admissible coframe with $d\theta=ih_{\alpha\overline{\beta}}\theta^\alpha\wedge\theta^{\overline{\beta}}$, then the \emph{connection forms} $\omega\indices{_\alpha^\beta}$ and the \emph{torsion forms} $\tau_\alpha =A_{\alpha\beta}\theta^\beta$ of the Tanaka-Webster connection are defined by $$d\theta^\beta=\theta^\alpha\wedge\omega\indices{_\alpha^\beta}+\theta\wedge\tau^\beta,\qquad \omega_{\alpha\overline{\beta}}+\omega_{\overline{\beta}\alpha}=dh_{\alpha\overline{\beta}},\qquad A_{\alpha\beta}=A_{\beta\alpha},$$ where indices are raised and lowered with $h_{\alpha\overline{\beta}}$, \emph{i.e.} $\omega_{\alpha\overline{\beta}}=h_{\sigma\overline{\beta}}\omega\indices{_\alpha^\sigma}$ \cite{Web78,Lee88}. We then have, for the dual frame  $(T_1,\ldots,T_n)$ to $(\theta^1,\ldots,\theta^n)$, $\nabla^\theta T_\alpha=\omega\indices{_\alpha^\beta}\otimes T_\beta$.

%This connection is called the \emph{Tanaka-Webster connection} of $(M,H,J,\theta)$. 
Due to the first condition in Theorem \ref{tw}, the torsion of the Tanaka-Webster connection is nonvanishing; however, we define:

\begin{defi}
The \emph{pseudohermitian torsion} $\tau$ of the Tanaka-Webster connection is the operator $\tau=\iota_R T^\theta$. If $\tau$ vanishes, $(M,H,J,\theta)$ is called \emph{normal}.% $\theta$ is then called a \emph{normal contact form}.
\end{defi}

\noindent Note that the definition of Tanaka-Webster connection implies that the pseudohermitian torsion is always trace-free as an endomorphism of the real vector bundle $H$.

%A strictly pseudoconvex pseudohermitian manifold whose pseudohermitian torsion vanishes is called \emph{normal}.

Let $\mathcal{R}^\theta$ be the curvature tensor field corresponding to the Tanaka-Webster connection. It can be decomposed into vertical, mixed, and horizontal terms. The vertical and mixed terms only depend on $\tau$ and its first derivatives. The horizontal part gives the \emph{Webster curvature tensor}. Let $\mathrm{Ric}_W(J,\theta)$ be its Ricci tensor, and $\mathrm{Scal}_W(J,\theta)$ be its scalar curvature, called the \emph{Webster scalar curvature}. In other words, the \emph{curvature forms} $\Pi\indices{_\alpha^\beta}=d\omega\indices{_\alpha^\beta}-\omega\indices{_\alpha^\sigma}\wedge\omega\indices{_\sigma^\beta}$ verify $$\Pi\indices{_\alpha^\beta}=\mathcal{R}\indices{^\theta_\alpha^\beta_\rho_{\overline{\sigma}}}\theta^\rho\wedge\theta^{\overline{\sigma}}\ \text{mod}\ \theta.$$ We then have $$\mathrm{Ric}_W(J,\theta)(T_\alpha,T_{\overline{\beta}})=\mathcal{R}\indices{^\theta_\alpha^\rho_\rho_{\overline{\beta}}}\quad \mathrm{and}\quad \mathrm{Scal}_W(J,\theta)=h^{\alpha\overline{\beta}}\mathrm{Ric}_W(J,\theta)(T_\alpha,T_{\overline{\beta}}).$$

\begin{defi}[\cite{CY13,Wan15}]
A strictly pseudoconvex pseudohermitian manifold $(M,H,J,\theta)$ is said to be \emph{pseudo-Einstein} if
\begin{displaymath}
%\left\lbrace
%\begin{array}{ll}
\mathrm{Ric}_W=\frac{1}{n}\mathrm{Scal}_W\gamma \quad\text{if}\ n\geq 2,\qquad
\mathrm{Scal}_{W,1}=i\tau_{1,\overline{1}}^{\overline{1}} \quad\text{if}\ n=1.
%\end{array}
%\right.
\end{displaymath}
% $\theta$ is then called a \emph{pseudo-Einstein contact form}.
%\end{defi}
%Case-Yang section 3
%\noindent In particular, if $(M,H,J)$ is embeddable, then it admits a pseudo-Einstein structure \cite{Lee88}. Note that, if $n\geq 2$ and $(M,H,J,\theta)$ is pseudo-Einstein, then the CR Bianchi identity implies that $$\mathrm{Scal}_{W,\alpha}=i(n-1)\tau_{\alpha,\overline{\beta}}^{\overline{\beta}}.$$
%Therefore, a pseudo-Einstein strictly pseudoconvex pseudohermitian manifold does not necessarily have constant scalar Webster curvature. However, a \emph{normal} pseudo-Einstein strictly pseudoconvex pseudohermitian manifold has constant scalar Webster curvature. A normal strictly pseudoconvex pseudohermitian $3$-manifold is pseudo-Einstein if and only if it has constant Webster scalar curvature.
%\begin{defi}[\cite{Wan15}]
A normal, pseudo-Einstein contact form is called an \emph{Einstein contact form}.
\end{defi}

\begin{ex}
The sphere $\Sph^{2n+1}\subset\C^{n+1}$ can be endowed with the contact form $$\theta_0=i\left(z_jd\overline{z}^j-\overline{z}_jdz^j\right)|_{\Sph^{2n+1}}.$$ The induced CR structure $(\Sph^{2n+1},H_0,J_0)$ is called the \emph{standard CR structure} of $\Sph^{2n+1}$. The pseudohermitian manifold $(\Sph^{2n+1},H_0,J_0,\theta_0)$ is Einstein, with constant positive Webster scalar curvature $\mathrm{Scal}_W(J_0,\theta_0)= \frac{n(n+1)}{2}$.  A CR manifold is called \emph{spherical} if it is locally CR equivalent to $(\Sph^{2n+1},H_0,J_0)$.
\end{ex}

\subsection{Circle bundles over a Riemann surface}

\label{circbun}

We recall here a construction detailed by D. Burns and C. Epstein \cite{BE88}, that will be useful in Section \ref{comp}. Let us consider a compact Riemann surface $\Sigma$ with a Hermitian metric $\gamma$. Let $T^{1,0}\Sigma$ be the holomorphic tangent bundle to $\Sigma$, and let $M$ be the unit circle bundle in $T^{1,0}\Sigma$. $M$ is then a $U(1)$-bundle over $\Sigma$, whose dual coframe gives a canonical one-form $\Theta^1$ on $M$. Moreover, since $\dim\Sigma=2$, $\gamma$ is automatically a Kähler metric, hence there is a unique torsion-free connection form $\Theta_1^1$ such that $d\Theta^1=\Theta^1\wedge \Theta^1_1$. We then have $$d\Theta^1_1=K\Theta^1\wedge \Theta^{\overline{1}},$$ where $K$ is the Gauss curvature of $\Sigma$. If $K$ never vanishes, an associated normal strictly pseudoconvex pseudohermitian structure $(J,\theta)$ is given on $M$ by $\theta=i\ \mathrm{sign}(K)\Theta_1^1$ and $\theta^1=\sqrt{|K|}\Theta^1$, so that $d\theta=i\theta^1\wedge \theta^{\overline{1}}$. Moreover, we have $$d\theta^1=\theta^1\wedge\left(\frac{1}{2}\frac{|K|_{,1}}{|K|}\theta^1-\frac{1}{2}\frac{|K|_{,\overline{1}}}{|K|}\theta^{\overline{1}}-i\ \mathrm{sign}(K)\theta\right).$$ If $K$ is constant, then $\omega_1^1=-i\ \mathrm{sign}(K)\theta$, hence $\mathrm{Scal}_W(J,\theta)=\mathrm{sign}(K)$.

Note that, by the following result, all non-spherical SPCR compact $3$-manifolds which admit a normal contact form are such bundles or finite quotients of them, \emph{i.e.} Seifert bundles.

\begin{prop}[\cite{Bel01}]
Let $(M,H,J)$ be a compact normal SPCR $3$-manifold. Then $(M,H,J)$ is either a finite quotient of the standard sphere or of a circle bundle over a Riemann surface of positive genus.
\end{prop}

\section{The contact Yamabe invariant}
\label{CRYam}

\subsection{The CR Yamabe problem}
\label{CRYam}

%\subsection{Presentation}

Let $(M,H,J,\theta)$ be a compact strictly pseudoconvex pseudohermitian manifold of dimension $2n+1$. % Let $\nabla$ be its Tanaka-Webster connection. 
We already mentioned that the set of positive contact forms on $(M,H,J)$ is a conformal class: $$[\theta]=\lbrace u^{\frac{2}{n}}\theta \ | \ u\in C^\infty(M,\R_+^*)\rbrace.$$ Here, the choice of the exponent $\frac{2}{n}$ is made to simplify further conformal change formulas.

\begin{defi}
Let $\pi_b:TM\rightarrow H$ be the pseudohermitian projection. The \emph{horizontal gradient} is the operator $\nabla_b:=\pi_b\nabla^\theta$. The \emph{sublaplacian} is $\Delta_b:=\mathrm{div}(\nabla_b\cdot)$.
\end{defi}

The similarity between conformal and CR geometry can be seen through the variation of the Webster scalar curvature under conformal changes of $\theta$: given a conformal factor $u$ in $C^\infty(M,\R_+^*)$, we have $$\mathrm{Scal}_W(J,u^{\frac{2}{n}}\theta)=u^{-\frac{n+2}{n}}\left(2\frac{n+1}{n}\Delta_b u +\mathrm{Scal}_W(J,\theta)u\right).$$ Therefore, $u^{\frac{2}{n}}\theta$ has constant Webster curvature $\mathrm{Scal}_W\equiv\lambda$ if and only if $$2\frac{n+1}{n}\Delta_b u +\mathrm{Scal}_W(J,\theta)u = \lambda u^{\frac{n+2}{n}},$$ which we will call the \emph{CR Yamabe equation}. This equation may be compared with the Riemannian Yamabe equation for a manifold of dimension $2n+2$. %This similarity supports the homogeneity mentioned in Section \ref{heis}.

%Let us denote $$P_\theta:=\Delta_H +\frac{n}{2(n+1)}\mathrm{Scal}_W(J,\theta)$$ the \emph{CR sublaplacian} on $C^\infty(M,\R)$. This operator is pseudo-conformally covariant: $$P_{u^{2/n}\theta}f=u^{-\frac{n+2}{n}}P_\theta (u) f \quad\forall f\in C^{\infty} (M,\R).$$
By analogy with the conformal case, the CR Yamabe problem is the following question: is there a constant Webster scalar curvature positive contact form in the conformal class $[\theta]$?

As in the conformal case, a sufficient condition is that there exists a contact form which realizes the infimum of the CR invariant 
\begin{equation}
\label{ycrdef}
Y_{CR}(M,H,J):=\inf_{\hat{\theta}\in [\theta]_1}S_W(M,J,\hat{\theta}),
\end{equation} 
where $$[\theta]_1:=\lbrace \hat{\theta}\in[\theta], \ \mathrm{Vol}(M,\hat{\theta}):=\int_M\hat{\theta}\wedge d\hat{\theta}^n=1\rbrace,$$ and $$S_W(M,J,\hat{\theta}):=\int_M\mathrm{Scal}_W(J,\hat{\theta})\hat{\theta}\wedge d\hat{\theta}^n$$ denotes the integral Webster scalar curvature. The functional $Y_{CR}$ is maximal for the standard sphere:

\begin{theo}[\cite{JL87}]\label{jl}
$Y_{CR}(M,H,J)\leq Y_{CR}(\Sph^{2n+1},H_0,J_0)=2\pi n(n+1)$.
\end{theo}

A positive contact form minimizing $Y_{CR}$ is called a \emph{Yamabe contact form}. The CR Yamabe problem has been given a positive answer by the following results of D. Jerison and J. Lee, and N. Gamara and R. Yacoub:

\begin{theo}[\cite{JL87}]\label{jlb}
If $Y_{CR}(M,H,J)<Y_{CR}(\Sph^{2n+1},H_0,J_0)$, then there is a Yamabe contact form.
\end{theo}

\begin{theo}[\cite{JL89}]
\label{jl2}
If $n\geq 2$ and $(M,H,J)$ is not spherical, then $Y_{CR}(M,H,J)<Y_{CR}(\Sph^{2n+1},H_0,J_0)$.
\end{theo}

\begin{theo}[\cite{Gam01,GY01}]
If $n=1$ or $(M,H,J)$ is spherical, then the CR Yamabe problem has a solution.
\end{theo}

\noindent The proof of this last theorem uses a technique of critical points at infinity initiated by A. Bahri. Note that the positive contact forms found this way are not necessarily Yamabe contact forms. However, J.-H. Cheng, A. Malchiodi and P. Yang have shown that Yamabe contact forms always exist on SPCR $3$-manifolds with non-negative CR Paneitz operator \cite{CMY17}. %, cf. Section \ref{crpan}.

On Einstein strictly pseudoconvex pseudohermitian manifolds, the following result by X. Wang ensures that all constant Webster scalar curvature contact forms are Einstein:

\begin{theo}[\cite{Wan15}]\label{einwang}
Let $(M,H,J)$ be a compact SPCR manifold which admits an Einstein contact form $\theta$. If $\tilde{\theta}=u^\frac{2}{n}\theta\in[\theta]$ has constant Webster scalar curvature, then $\tilde{\theta}$ is Einstein. Moreover, if $(M,H,J)$ is non-spherical, then $u$ is constant.
\end{theo}

\subsection{The contact Yamabe invariant}

The resolution of the CR Yamabe problem, cf. Section \ref{CRYam}, leads naturally to the consideration of the following quantity:

\begin{defi}[\cite{Wu09}]
\label{sigmaCR}
Let $(M,H)$ be a compact SPC manifold. Let $\mathcal{J}$ be the set of complex structures $J$ on $(M,H)$ such that $(M,H,J)$ is SPCR. The \emph{contact Yamabe invariant} $\sigma_c(M,H)$ is defined by $$\sigma_c(M,H):=\sup_{\mathcal{J}}Y_{CR}(M,H,J).$$
\end{defi}
%, defined for an SPC manifold $(M,H)$ by $$\sigma_c(M,H):=\sup_{\mathcal{J}} \inf_{\hat{\theta}\in [\theta]_1}\int_M\mathrm{Scal}_W(J,\hat{\theta})\ \hat{\theta}\wedge (d\hat{\theta})^n,$$ where $\mathrm{Scal}_W$ denotes the Webster scalar curvature, the supremum runs over the set $\mathcal{J}$ of all complex structures $J$ on $(M,H)$ such that $(M,H,J)$ is strictly pseudoconvex and the infimum runs over all compatible contact forms of volume 1 on $(M,H)$ \cite{Wu09}.
 As mentioned in the introduction, $\sigma_c$ is %a global contact invariant of $(M,H)$. %It has been introduced by C.-T. Wu, but to our knowledge it has not been used afterwards \cite{Wu09}. This invariant is 
an actual contact invariant in dimension $3$, in the sense that, since $\mathcal{J}$ is always non-empty, all contact $3$-manifolds are SPC. Note that few contact invariants are currently available: they are necessarily global by Darboux's theorem, and most of them come from homological considerations. In higher dimension, for some contact structures, due to the obstructions on the integrability of complex structures and on their compatibility with a given contact form, the set $\mathcal{J}$ might be empty.

%In the conformal setting, the Yamabe invariant, or Yamabe invariant, has been introduced by R. Schoen \cite{Sch89} and O. Kobayashi \cite{Kob87}. In the CR setting, a corresponding invariant has been defined by C.-T. Wu.

%We now present a rigorous definition of $\sigma_c$. The definition (\ref{ycrdef}) of $Y_{CR}$ can be extended to noncompact manifolds via the following result:

As in the conformal case, the contact Yamabe invariant characterizes manifolds which admit a structure with positive curvature:

\begin{prop}[\cite{Wan03}]\label{pos}
%Let $(M,H,J)$ be an SPCR manifold. Then $Y_{CR}(M,H,J)>0$ if and only if there exists a positive contact form $\theta$ on $(M,H,J)$ with positive Webster scalar curvature.
Let $(M,H)$ be a compact SPC manifold. Then $\sigma_c(M,H)>0$ if and only if there exists a strictly pseudoconvex pseudohermitian structure $(J,\theta)$ on $(M,H)$ with positive Webster scalar curvature.
\end{prop}
%\begin{proof}
%Cf prop 5.1 in \cite{Wan03}.
%\end{proof}

Finally, let us recall the following lemma, that will be essential in Section \ref{glu}.

\begin{lemm}
\label{lip}
Let $(M,H,J)$ be an SPCR manifold. The infimum in Definition (\ref{ycrdef}) of $Y_{CR}(M,H,J)$ may be taken over the space $L_c(M,\R_+)$ of non-negative Lipschitz functions with compact support on $M$.
\end{lemm}

\begin{proof} %(\cite{JL87} lemma 4.5, \cite{LP87} pp.49-50).\\
Indeed, $Y_{CR}(M,H,J)=\underset{u\in C^\infty(M,\R_+^*)}{\inf}Q_\theta(u)$ where $\theta$ is any positive contact form on $(M,H,J)$ and $$Q_\theta(u):=\frac{\int_M \left(2\frac{n+1}{n}|\nabla u|^2+\mathrm{Scal}_W(J,\theta) u^2\right) \theta\wedge d\theta^n}{\left(\int_M u^{2\frac{n+1}{n}} \theta\wedge d\theta^n\right)^\frac{n}{n+1}},$$ hence $Q_\theta$ is continuous in the Sobolev space $W^{1,2}(M)$. Since $C^\infty(M)$ is dense in $W^{1,2}(M)$, since $Q_\theta(|u|)=Q_\theta(u)$ for all $u\in C^\infty(M)$, and since a nonnegative Lipschitz fonction can be arbitrarily approximated in $W^{1,2}$ norm by a positive smooth function, we have $Y_{CR}(M,H,J)=\underset{u\in L_c(M,\R_+)}{\inf}Q_\theta(u)$.
\end{proof}

\section{CR handle attaching on a spherical manifold}
\label{glu2}

We recall here a handle attaching process on spherical SPCR manifolds, compatible with the CR structure, which is due to W. Wang \cite{Wan03}. If the handle is attached between two distinct connected components, this provides a connected sum of the components.%, which we call a \emph{Wang connected sum}. Otherwise, we call it a \emph{Wang handle attaching}.

Let either $\left(M=M_1\sqcup M_2,H,J,\hat{\theta}\right)$ be a disjoint union of two connected spherical strictly pseudoconvex pseudohermitian manifolds, and $p_1\in M_1$, $p_2\in M_2$; or $\left(M,H,J,\hat{\theta}\right)$ be a connected spherical strictly pseudoconvex pseudohermitian manifold, and $p_1, p_2\in M$.

Let $M_0=M\setminus\lbrace p_1,p_2\rbrace$. For $i$ in $\lbrace 1,2\rbrace$, let $U_i\subset M$ be a neighbourhood of $p_i$ and $$\varphi_i:U_i\rightarrow B(0,2):=\lbrace \xi\in\Heis^{2n+1}, \ \|\xi\|_\Heis<2\rbrace$$ local coordinates such that $\varphi_i(p_i)=0.$ Let us denote, for $0<r<1$, $$U_i(r)=\lbrace x\in U_i, \ \|\varphi_i(x)\|_\Heis<r\rbrace,$$ $$U_i(r,1)=\lbrace x\in U_i, \ r<\|\varphi_i(x)\|_\Heis<1\rbrace.$$

Since $M$ is spherical around $p_1$ and $p_2$, there exists $\lambda\in C^\infty(M_0,\R_+^*)$ such that, denoting $\theta=\lambda\hat{\theta}$ on $M_0$, $$\forall i\in\lbrace 1,2\rbrace,\ \forall \xi\in B(0,2)\setminus\lbrace 0\rbrace, \quad {\varphi_i}_*\theta(\xi)=\|\xi\|_\Heis^{-2}\theta_{\Heis}(\xi),$$ that is, $(M_0,H,J,\theta)$ has cylindrical ends. Indeed, we can define a mapping 
\begin{displaymath}
\begin{array}{cccc}
\Phi: & B(0,1) & \longrightarrow & \R_+\times\Sigma^{2n}\\
& \xi & \mapsto & \left(\log\frac{1}{\xi},\frac{\xi}{\|\xi\|_\Heis}\right),
\end{array}
\end{displaymath}
where $\Sigma^{2n}:=\lbrace\xi\in\Heis^{2n+1}, \ \|\xi\|_\Heis=1\rbrace$, and $\tilde{\theta}:=\Phi_*(\|\cdot\|_\Heis^{-2}\theta_{\Heis})$. Then $$(B(0,1),H_{\Heis},J_{\Heis},\|\cdot\|_\Heis^{-2}\theta_{\Heis})\simeq(\R_+\times\Sigma^{2n},\tilde{H},\tilde{J},\tilde{\theta}),$$ where the equivalence is pseudohermitian. Denoting $\hat{M}:=M\setminus U_1(1)\cup U_2(1)$, 
\begin{equation}
(M_0,H,J,\theta)\simeq(\R_+\times\Sigma^{2n},\tilde{H},\tilde{J},\tilde{\theta})\ \cup\ (\hat{M},H,J,\theta)\ \cup\ (\R_+\times\Sigma^{2n},\tilde{H},\tilde{J},\tilde{\theta}).
\end{equation}

Now, let us denote, for $r\in(0,1)$ and $A\in U(n)$, by $\psi_{r,A}:U_1(r,1)\rightarrow U_2(r,1)$ the mapping $$\psi_{r,A}=\varphi_2^{-1}\circ\delta_r\circ R\circ U_A\circ\varphi_1,$$ where $$\delta_r:(z,t)\mapsto(rz,r^2 t),$$ $$U_A:(z,t)\mapsto(Az,t),$$ and $$R:(z,t)\mapsto\left(\frac{-z}{|z|^2-it},\frac{-t}{|z|^4+t^2}\right)$$ denote respectively dilations, unitary transformations and inversion in $\Heis^{2n+1}$.

\noindent Let $(M_{r,A},H_{r,A},J_{r,A},\theta_{r,A})$ be the pseudohermitian manifold formed from $M$ by removing $\overline{U_1(r)}$ and $\overline{U_2(r)}$, and by identifying $U_1(r,1)$ with $U_2(r,1)$ along $\psi_{r,A}$. Let $$\pi_{r,A}:M\setminus\overline{U_1(r)\cup U_2(r)}\rightarrow M_{r,A}$$ be the corresponding projection.

\noindent Since $$\delta_r^*{\varphi_i}_*\theta=\frac{1}{r^2\|\cdot\|_\Heis^2}r^2\theta_{\Heis}={\varphi_i}_*\theta$$ and $$R^*{\varphi_i}_*\theta=\frac{1}{\|R(\cdot)\|_\Heis^2}R^*\theta_{\Heis}=\frac{\|\cdot\|_\Heis^2}{\|\cdot\|_\Heis^4}\theta_{\Heis}={\varphi_i}_*\theta,$$ the gluing preserves $\theta$ on $U_i(r,1)$. Hence, $$\pi_{r,A}^*\theta_{r,A}=\theta \qquad \mathrm{on} \ M\setminus\overline{U_1(r)\cup U_2(r)}.$$ We have in fact 
\begin{equation}
\label{dec1}
(M_{r,A},H_{r,A},J_{r,A},\theta_{r,A})\simeq(\hat{M},H,J,\theta)\ \cup\ ([0,l]\times\Sigma^{2n},\tilde{H},\tilde{J},\tilde{\theta})
\end{equation} where $l=\log\left(\frac{1}{r}\right)\in(0,+\infty)$.

\section{A CR Kobayashi inequality}
\label{CRadh}

\subsection{Non-decreasing of $\sigma_c$ under SPC handle attaching}
\label{glu}

This part follows a method developed in the conformal setting by O. Kobayashi \cite{Kob87}. It has been adapted in the CR setting by W. Wang, and more recently by J.-H. Cheng and H.-L. Chiu \cite{Wan03,CC19}. A similar technique has been implemented in the quaternionic context \cite{SW16}. Theorem \ref{CRkobh} is a direct consequence of the following result:% and of Proposition \ref{eqlem}:

\begin{theo}\label{YCRkob}
Let $(M,H,J)$ be a compact SPCR manifold. Let $(\tilde{M},\tilde{H},\tilde{J})$ be a manifold obtained from $(M,H,J)$ by CR handle attaching, then $$Y_{CR}(\tilde{M},\tilde{H},\tilde{J})\geq Y_{CR}(M,H,J).$$
%Let $(M_1,H_1,J_1)$ and $(M_2,H_2,J_2)$ be two compact SPCR manifolds of dimension $2n+1$. Let $(M_1,H_1,J_1)\# (M_2,H_2,J_2)$ be their connected sum, then $$Y_{CR}\left((M_1,H_1,J_1)\# (M_2,H_2,J_2)\right)\geq Y_{CR}\left((M_1,H_1,J_1)\sqcup (M_2,H_2,J_2)\right).$$
\end{theo}

\begin{proof}
Let $p_1,p_2\in M$. We use the following lemma, that will be proved in Section \ref{lsph}:

\begin{lemm}
\label{locsphe}
We may assume that $(M,H,J)$ is spherical around $p_1$ and $p_2$.
\end{lemm}
\noindent Under this assumption, we can apply the construction of Section \ref{glu2} to $M$. Let $(M_{r,A},H_{r,A},J_{r,A})$ be obtained from $(M,H,J)$ by CR handle attaching.

By definition of $Y_{CR}(M_{r,A},H_{r,A},J_{r,A})$, there exists a function $f_l\in C^\infty(M_{r,A},\R_+^*)$ such that
\begin{align}
\displaystyle{S_W\left(M_{r,A},J_{r,A},f_l^{2/n}\theta_{r,A}\right)} & =\displaystyle{\int_{M_{r,A}}\left(2\left(\frac{n+1}{n}\right)|\nabla^{\theta_{r,A}} f_l|^2+\mathrm{Scal}_W(J_{r,A},\theta_{r,A})f_l^2\right) \theta_{r,A}\wedge d\theta_{r,A}^n} \notag\\
& \displaystyle{<Y_{CR}(M_{r,A},H_{r,A},J_{r,A})+\frac{1}{1+l}} \label{yamdef}
\end{align}
and $$\int_{M_{r,A}}f_l^{2\left(\frac{n+1}{n}\right)}\theta_{r,A}\wedge d\theta_{r,A}^n=1.$$

\begin{lemm}
\label{majlem}
There exists $l_*\in[0,l]$ such that $$\int_{\lbrace l_*\rbrace\times\Sigma^{2n}}\left(|d_b f_l|^2+f_l^2\right) dS_{\Sigma^{2n}}\leq\frac{C}{l},$$ where $C$ is a constant independent of $l$.
\end{lemm}

\begin{proof}
Let $C_1=-\min(0,\underset{\hat{M}}{\min}\ \mathrm{Scal}_W(J,\theta))\mathrm{Vol}(\hat{M},\theta)^\frac{1}{n+1}$. H\"{o}lder's inequality yields $$\int_{\hat{M}}f_l^2\theta\wedge d\theta^n\leq\mathrm{Vol}(\hat{M},\theta)^\frac{1}{n+1},$$ and then, using decomposition (\ref{dec1}), $$\int_{[0,l]\times\Sigma^{2n}}\left(2\left(\frac{n+1}{n}\right)|d_b f_l|^2+\mathrm{Scal}_W(\tilde{J},\tilde{\theta})f_l^2\right) \tilde{\theta}\wedge d\tilde{\theta}^n\leq Y_{CR}(M_{r,A},H_{r,A},J_{r,A})+\frac{1}{1+l}+C_1.$$ Consequently there exists $l_*\in[0,l]$ such that $$\int_{\lbrace l_*\rbrace\times\Sigma^{2n}}\left(2\left(\frac{n+1}{n}\right)|d_b f_l|^2+\mathrm{Scal}_W(\tilde{J},\tilde{\theta})f_l^2\right) \tilde{\theta}\wedge d\tilde{\theta}^n\leq\frac{1}{l}\left(Y_{CR}(M_{r,A},H_{r,A},J_{r,A})+\frac{1}{1+l}+C_1\right).$$ The lemma is obtained with $C=\displaystyle{\frac{Y_{CR}(M_{r,A},H_{r,A},J_{r,A})+1+C_1}{\min\left(2\left(\frac{n+1}{n}\right),\underset{\lbrace l_*\rbrace\times\Sigma^{2n}}{\min}\mathrm{Scal}_W(\tilde{J},\tilde{\theta})\right)}}$.
\end{proof}

We therefore decompose
\begin{equation}
\begin{array}{rcl}
(M_0,H,J,\theta)&\simeq&([l_*,\infty)\times\Sigma^{2n},\tilde{H},\tilde{J},\tilde{\theta})\\
&& \cup\ (M_{r,A}\setminus \lbrace l_*\rbrace\times \Sigma^{2n},H_{r,A},J_{r,A},\theta_{r,A})\\
&& \cup\ ([l-l_*,\infty)\times\Sigma^{2n},\tilde{H},\tilde{J},\tilde{\theta})
\end{array}
\end{equation}
and extend $f_l$ to $M_0$ as follows: $F_l=f_l$ on $M_{r,A}\setminus \lbrace l_*\rbrace\times \Sigma^{2n}$ and
\begin{displaymath}
F_l:(s,\xi)\mapsto\left\lbrace
\begin{array}{ll}
(l_*+1-s)f_l(l_*,\xi) & \forall (s,\xi)\in[l_*,l_*+1]\times\Sigma^{2n},\\
0 &  \forall (s,\xi)\in[l_*+1,\infty)\times\Sigma^{2n},\\
0 &  \forall (s,\xi)\in[l-l_*+1,\infty)\times\Sigma^{2n},\\
(l-l_*+1-s)f_l(l_*,\xi) &  \forall (s,\xi)\in[l-l_*,l-l_*+1]\times\Sigma^{2n}.\\
\end{array}\right.
\end{displaymath}

We thus obtain from (\ref{yamdef}) and Lemma \ref{majlem}
\begin{displaymath}
\begin{array}{llll}
\displaystyle{S_W(M_0,J,F_l^{2/n}\theta)} & = & \displaystyle{S_W\left(M_{r,A},J_{r,A},f_l^{2/n}\theta_{r,A}\right)}\\
& & + \displaystyle{\int_{[l_*,l_*+1]\times\Sigma^{2n}} \left(2\left(\frac{n+1}{n}\right)|d_b F_l|^2+\mathrm{Scal}_W(\tilde{J},\tilde{\theta})F_l^2\right) \tilde{\theta}\wedge d\tilde{\theta}^n}\\
& & + \displaystyle{\int_{[l-l_*,l-l_*+1]\times\Sigma^{2n}} \left(2\left(\frac{n+1}{n}\right)|d_b F_l|^2+\mathrm{Scal}_W(\tilde{J},\tilde{\theta})F_l^2\right) \tilde{\theta}\wedge d\tilde{\theta}^n}\\
& = & \displaystyle{S_W\left(M_{r,A},J_{r,A},f_l^{2/n}\theta_{r,A}\right)}\\
& & + 2 \displaystyle{\int_{\Sigma^{2n}} \left(\frac{2}{3}\frac{n+1}{n}|d_b f_l(l_*,\cdot)|^2 + \left(2\left(\frac{n+1}{n}\right)+\frac{1}{3}\mathrm{Scal}_W(\tilde{J},\tilde{\theta})\right)f_l(l_*,\cdot)^2 \right) dS_{\Sigma^{2n}}}\\
& \leq & \displaystyle{S_W\left(M_{r,A},J_{r,A},f_l^{2/n}\theta_{r,A}\right)}\\
& & + 2 \displaystyle{\int_{\lbrace l_*\rbrace\times\Sigma^{2n}} \left(\frac{2}{3}\frac{n+1}{n}|d_b f_l|^2 + \left(2\left(\frac{n+1}{n}\right)+\frac{1}{3}\mathrm{Scal}_W(\tilde{J},\tilde{\theta})\right)f_l^2 \right) dS_{\Sigma^{2n}}}\\
& \leq & \displaystyle{Y_{CR}(M_{r,A},H_{r,A},J_{r,A})+\frac{B}{l}},
\end{array}
\end{displaymath}
where $B$ is a constant independent of $l$, and $$\int_{M_0}F_l^{2\left(\frac{n+1}{n}\right)}\theta\wedge d\theta^n>1.$$

Since the infimum in the Yamabe functional may be taken over all nonnegative Lipschitz functions with compact support as conformal factors %(\cite{JL87} lemma 4.5, \cite{LP87} pp.49-50)
by Lemma \ref{lip}, we get that $$Y_{CR}(M_0,H,J)\leq Y_{CR}(M_{r,A},H_{r,A},J_{r,A})+\frac{B}{l},$$ which, for $l$ sufficiently large, yields the desired inequality.
\end{proof}

%Replacing $M_1\sqcup M_2$ with a connected manifold $M$ and taking $p_1, p_2 \in M$ in the last proof, we can identically show the following result:
%\begin{prop}
%The CR Yamabe invariant is nondecreasing under handle attaching.% on an SPC manifold $(M,H)$.
%\end{prop}

\subsection{Local sphericity}
\label{lsph}

In this section, we prove the following technical lemma, which is essential for the proof of Theorem \ref{YCRkob}.
%Following \cite{Kob87}, we prove that:

\begin{lemm}
\label{locsphe}
%We may assume in Theorem \ref{YCRkob} that $(M_1,H_1,J_1)$ and $(M_2,H_2,J_2)$ are spherical around the connecting points. Namely, g
Let $(M,H,J)$ be a compact SPCR manifold. Given $p_1$ and $p_2$ in $M$, there is a 1-parameter family of complex structures $(J_t)$ in $\mathcal{J}$ $C^0$-converging to $J$ such that $(M,H,J_{t})$ is spherical around $p_1$ and $p_2$, and $Y_{CR}(M,H,J_{t})\underset{t\rightarrow 0}{\longrightarrow}Y_{CR}(M,H,J)$.
%with connected sum $(M_{r_t,A_t},H_{r_t,A_t},J_{r_t,A_t})$
\end{lemm}

\noindent In other words, to prove Theorem \ref{YCRkob}, we may assume that $(M,H,J)$ is spherical around $p_1$ and $p_2$. This lemma is a direct consequence of the two following results:

\begin{lemm}
\label{cont}
Let $(M,H,J,\theta)$ be a compact strictly pseudoconvex pseudohermitian manifold. Let $(J_t)$ be a $1$-parameter family of complex structures in $\mathcal{J}$ $C^0$-converging to $J$ such that $\mathrm{Scal}_W(J_t,\theta)\underset{t\rightarrow 0}{\longrightarrow}\mathrm{Scal}_W(J,\theta)$. Then $Y_{CR}(M,H,J_t)\underset{t\rightarrow 0}{\longrightarrow}Y_{CR}(M,H,J)$.
\end{lemm}

\begin{lemm}[\cite{CCH19}]
\label{locsphe2}
Let $(M,H,J,\theta)$ be a compact strictly pseudoconvex pseudohermitian manifold. Let $p_1$ and $p_2$ be two points in $M$. Let  $(\varepsilon_t)$ be a 1-parameter family of positive numbers decreasing to $0$. There is a 1-parameter family of complex structures $(J_t)$ in $\mathcal{J}$ $C^0$-converging to $J$ such that for all $t$, $J_t$ coincides with $J$ outside an $\varepsilon_t$-neighbourhood $U'_t$ of $\lbrace p_1,p_2\rbrace$, $J_t$ is spherical inside $U_t\subset U'_t$, and $|\mathrm{Scal}_W(J_t,\theta)-\mathrm{Scal}_W(J,\theta)|\leq \varepsilon_t$.
\end{lemm}
%The proof of Lemma \ref{locsphe2} is due to H.-L. Chiu \cite{Chi}.

\begin{proof}[Proof of Lemma \ref{cont}]
We adapt from the conformal case a proof due to L. B\'{e}rard Bergery \cite{Ber83}. 
%Since $Y_{CR}(M,H,\cdot)$ is the infimum of a continuous functional on $\mathcal{J}$, it is upper semi-continuous. Let us show that it is also lower semi-continuous.
%Let $\theta$ be a unit volume contact form on $(M,H)$ and let $J_t$ be a $C^2$-1-parameter family of complex structures in $\mathcal{J}$ such that $(J_t,\theta)$ is a strictly pseudoconvex pseudohermitian structure for all $t$.
Let us denote $$F=\left\lbrace u\in C^\infty(M,\R_+^*)\ \left|\ \int_M u^{2\left(\frac{n+1}{n}\right)}\theta\wedge d\theta^n=1\right.\right\rbrace.$$
By definition, $Y_{CR}(M,H,J_t)=\underset{u\in F}{\inf}I_t(u)$ where $$I_t(u)=2\left(\frac{n+1}{n}\right)\int_M |d_b u|_t^2 \theta\wedge d\theta^n+\int_M u^2\mathrm{Scal}_W(J_t,\theta) \theta\wedge d\theta^n.$$
Let $\varepsilon>0$. For each $t$ there exists $u_t$ in $F$ such that $$Y_{CR}(M,H,J_t)\leq I_t(u_t)\leq Y_{CR}(M,H,J_t)+\varepsilon.$$
%Moreover, since (e.g. \cite{CL90}, 2.20), %s 
%\begin{equation}
%\label{var1}\mathrm{Scal}_W(J,(1+\varepsilon)^{\frac{2}{n}}\theta_1)-\mathrm{Scal}_W(J,\theta_1)=2\left(\frac{n+1}{n}\right)\Delta_H\varepsilon-\frac{2}{n}\mathrm{Scal}_W(\theta_1)\varepsilon+o(\varepsilon)
%\end{equation}
%for $E\in\mathcal{E}(J):=\lbrace E\in \mathrm{End}(H) \ |\ EJ+JE=0\rbrace$, 
%$$\partial_J \mathrm{Scal}_W(J,\theta)(E)=\frac{i}{2}(E\indices{_\alpha^{\overline{\beta}}_{,\overline{\beta}}^\alpha}-E\indices{_{\overline{\alpha}}^\beta_{,\beta}^{\overline{\alpha}}})-\frac{1}{2}(\tau\indices{_\alpha^{\overline{\beta}}}
%E\indices{_{\overline{\beta}}^\alpha}+\tau\indices{_{\overline{\alpha}}^\beta}E\indices{_\beta^{\overline{\alpha}}}),$$ we have $\mathrm{Scal}_W(J_t,\theta)\underset{t\rightarrow 0}{\longrightarrow}\mathrm{Scal}_W(J_0,\theta)$. 
% where $\tau=\iota_T T^\nabla$ is the pseudohermitian torsion, 
%we finally get $$|\mathrm{Scal_W}(J_i,\theta)-\mathrm{Scal_W}(J,\theta)|\leq C_0\|\tau(J,\theta)\|\varepsilon_i^\frac{1}{2}+C_2\varepsilon_i+o(\varepsilon_i)\quad \text{on}\ M.$$
Let $\eta>0$ and $K>0$ be such that, for all $t$ in $[-\eta,\eta]$, $$|\mathrm{Scal}_W(J_t,\theta)-\mathrm{Scal}_W(J,\theta)|\leq\varepsilon,$$ $$|\mathrm{Scal}_W(J_t,\theta)|\leq K,$$ and $$\frac{1}{1+\varepsilon}\leq\frac{|\omega|^2_t}{|\omega|^2_0}\leq1+\varepsilon \quad \forall \omega \in\Omega^1(M).$$ In particular, $I_t(u_t)\leq I_t(1)+\varepsilon\leq K+\varepsilon$. Now, by H\"{o}lder's inequality, $$\int_M u_t^2\theta\wedge d\theta^n\leq 1,$$ hence we have
%\begin{displaymath}
%\begin{array}{ll}
%|I_t(u_t)-I_0(u_t)| & \leq 2\varepsilon\frac{n+1}{n}\int_M \|d u_t\|_t^2 \theta\wedge d\theta^n+\varepsilon\int_M u_t^2 \theta\wedge d\theta^n\\
% & \leq \varepsilon\left(I_t(u_t)-\int_M \mathrm{Scal}_W(J_t,\theta)u_t^2 \theta\wedge d\theta^n+1\right)\\
% & \leq \varepsilon(2K+\varepsilon+1).
%\end{array}
%\end{displaymath}
\begin{displaymath}
\begin{array}{ll}
I_0(u_t) & = 2\left(\frac{n+1}{n}\right)\int_M |d_b u_t|_0^2 \theta\wedge d\theta^n+\int_M u_t^2  \mathrm{Scal}_W(J_0,\theta)\theta\wedge d\theta^n\\
& \leq 2\left(\frac{n+1}{n}\right)(1+\varepsilon)\int_M |d_b u_t|_t^2 \theta\wedge d\theta^n+\int_M u_t^2 \mathrm{Scal}_W(J_t,\theta) \theta\wedge d\theta^n\\
& \quad +\int_M u_t^2|\mathrm{Scal}_W(J,\theta)-\mathrm{Scal}_W(J_t,\theta)|  \theta\wedge d\theta^n\\
& \leq 2\left(\frac{n+1}{n}\right)(1+\varepsilon)\int_M |d_b u_t|_t^2 \theta\wedge d\theta^n+\int_M u_t^2 \mathrm{Scal}_W(J_t,\theta) \theta\wedge d\theta^n+\varepsilon\\
 & \leq (1+\varepsilon)I_t(u_t)+\varepsilon(K+1),
\end{array}
\end{displaymath}
and similarly $$I_t(u_0)\leq (1+\varepsilon)I_0(u_0)+\varepsilon(K+1).$$

\noindent Then, for all $t$ in $[-\eta,\eta]$,
\begin{displaymath}
\begin{array}{ll}
\frac{1}{1+\varepsilon}Y_{CR}(M,H,J)-\frac{\varepsilon}{1+\varepsilon}(K+1)-\varepsilon & \leq Y_{CR}(M,H,J_t)\\
&\leq (1+\varepsilon)(Y_{CR}(M,H,J)+\varepsilon)+\varepsilon(K+1).
\end{array}
\end{displaymath}
\end{proof}

%\begin{lemm}
%\label{locsphe}
%Let $(M,H)$ be a compact SPC manifold and $\iota:\Sph^k\hookrightarrow M$ be an isotropic embedding. Let us call $S=\iota(\Sph^k)$ the \emph{surgery locus}. There exist a sequence of complex structures $(J_i)$ on $(M,H)$ $C^2$-converging to $J$ in $\mathcal{J}$ and a sequence of positive numbers $(\varepsilon_i)$ decreasing to $0$ such that for all $i$, $J_i$ coincides with $J_\mathrm{sph}$ spherical in an $\varepsilon_i$-neighbourhood $U_i$ of $S$, $J_i$ coincides with $J$ outside $U'_i\supset U_i$, and $\|\mathrm{Scal_W}(J_i,\theta)-\mathrm{Scal_W}(J,\theta)\|\leq \varepsilon_i$ on $M$.
%\end{lemm}

\begin{rk}
Since $Y_{CR}(M,H,J)$ only depends on derivatives up to order $2$ of $J$, the supremum in $\sigma_c(M,H)$ may be taken over all $C^2$ complex structures on $(M,H)$. Therefore, in the following proof, gluing complex structures only needs to be considered up to $C^2$-regularity.

\end{rk}

\begin{proof}[Proof of Lemma \ref{locsphe2} \cite{CCH19}]

%/!\ We suppose that $k=n$. If $k<n$, we must replace $C^c_\lambda$ by $C^a_\lambda$ in the proof to get submanifolds of the same dimension.

%Let $(M,H,J,\theta)$ be a strictly pseudoconvex pseudohermitian manifold and $\iota:\Sph^k\rightarrow M$ be an isotropic embedding with trivial conformal symplectic normal bundle. Let $\lambda$ be large enough so that $\mathrm{Scal}_W(J_\lambda,\theta_\lambda)$ is positive on $C^c_\lambda$. Let $\varepsilon>0$ and let $U$ be an $\varepsilon$-neighbourhood of $\iota(\Sph^k)$ in $M$ such that there is a strict contactomorphism $\Phi:U\rightarrow C^c_\lambda$, i.e. $\Phi^* \theta^c_\lambda=\theta$ ($\Phi$ exists by Weinstein's tubular neighbourhood theorem, \cite{Gei08} 2.5.8.; \cite{KM97}, 43.18.). %, denoting $U^c_\lambda=\Phi(U)$.
%We choose $\Phi$ such that $\Phi^* J^c_\lambda=J$ on $\iota(\Sph^k)$ (cf \cite{BR09}, 3.3).

We follow a construction due to O. Biquard and Y. Rollin \cite{BR09}. We assume that for all $t$, $B(p_1,\varepsilon_t)\cap B(p_2,\varepsilon_t)=\emptyset$, where the distances are taken with respect to the Webster metric. For a given $t$, let $U'_t$ be an $\varepsilon_t$-neighbourhood of $\lbrace p_1,p_2\rbrace$, and let $x=\min\left(d(\cdot,p_1),d(\cdot,p_2)\right)$ on $M$. There is a smooth cut-off function $w_t:\R_+\rightarrow\R_+$ such that $\chi_t:=w_t\circ x=0$ on some $U_t\subset U'_t$, $\chi_t=1$ outside $U'_t$, and for all $x$ in $\R_+$, $|xw'_t(x)|\leq\varepsilon_t$ and $|x^2 w''_t(x)|\leq\varepsilon_t$ (cf. \cite{Kob87}, Sublemma 3.4.). Indeed, we may take $w_t$ as a smoothing of $\tilde{w}_t$ defined by 
\begin{displaymath}
\tilde{w}_t(x)=\left\lbrace
\begin{array}{ll}
0 & \forall x\leq \varepsilon_t e^{-\frac{2}{\varepsilon_t}}\\
1-\frac{\varepsilon_t }{2}\log\left(\frac{\varepsilon_t}{x}\right) & \forall x\in [\varepsilon_t e^{-\frac{2}{\varepsilon_t}},\varepsilon_t]\\
1 & \forall x\geq \varepsilon_t.
\end{array}\right.
\end{displaymath}

If $\dim M=3$, then all almost complex structures are formally integrable. Let us take $i$ in $\lbrace 1,2\rbrace$. Let $\psi_i:U_i\rightarrow U_\Heis$ be a contactomorphism identifying a neighbourhood $U_i$ of $p_i$ in $M$ with a neighbourhood $U_\Heis$ of $0$ in $\Heis^3$ such that $\psi_i(p_i)=0$ and, denoting $\tilde{J}_i:=({\psi_i})_* J$ and $\tilde{\theta}_i:=({\psi_i})_* \theta$, such that $j^1_0(\tilde{J}_i)=j^1_0(J_\Heis)$ and $\mathrm{Scal}_W(\tilde{J}_i,\tilde{\theta}_i)=\mathrm{Scal}_W(J_\Heis,\tilde{\theta}_i)$ at $0$. We assume that $U_1\cap U_2=\emptyset$. For $t$ large enough that $U'_t\subset U_1\sqcup U_2$, let $\tilde{\chi}_{i,t}:=({\psi_i}^{-1})^*\chi_t|_{U_i}$. For $(z_1,y)$ in $\Heis^3$, let $\tilde{J}_{i,t}(z_1,y):=\tilde{J}_i(\tilde{\chi}_{i,t} z_1,\tilde{\chi}_{i,t}^2 y)$. Then $\tilde{J}_{i,t}$ coincides with $J_\Heis$ inside $\psi_i(U_t\cap U_i)$, and with $\tilde{J}_i$ outside $\psi_i(U'_t\cap U_i)$. Therefore, the complex structure $J_t$ defined on $M$ by $$\forall i\in\lbrace 1,2\rbrace,\ J_t:=\psi_i^*\tilde{J}_{i,t}\ \text{on}\ U_i, \quad J_t:=J\ \text{elsewhere},$$ has the desired properties. In particular, we have $|J_t-J|=O(x^2)$, $|\nabla(J_t-J)|=O(x)$, and $|\nabla^2(J_t-J)|=O(1)$. We then use Formula (4.7) in \cite{CT00}: %s 
%\begin{equation}
%\label{var1}\mathrm{Scal}_W(J,(1+\varepsilon)^{\frac{2}{n}}\theta_1)-\mathrm{Scal}_W(J,\theta_1)=2\left(\frac{n+1}{n}\right)\Delta_H\varepsilon-\frac{2}{n}\mathrm{Scal}_W(\theta_1)\varepsilon+o(\varepsilon)
%\end{equation}
for $E\in\mathcal{E}(J):=\lbrace E\in \mathrm{End}(H) \ |\ EJ+JE=0\rbrace$, we have $$\partial_J \mathrm{Scal}_W(J,\theta)(E)=\frac{i}{2}(E\indices{_{11,\overline{1}\overline{1}}}-E\indices{_{\overline{1}\overline{1},11}})-\frac{1}{2}(\tau\indices{_{11}} E\indices{_{\overline{1}\overline{1}}}+\tau\indices{_{\overline{1}\overline{1}}}E\indices{_{11}}).$$ 
In our case, for some constant $C$, and for $t$ sufficiently small, we then have %$\mathrm{Scal}_W(J_t,\theta)\underset{t\rightarrow 0}{\longrightarrow}\mathrm{Scal}_W(J_0,\theta)$. 
% where $\tau=\iota_T T^\nabla$ is the pseudohermitian torsion, 
\begin{displaymath}
\begin{array}{rcl}
|\mathrm{Scal}_W(J_t,\theta)-\mathrm{Scal}_W(J,\theta)|&\leq& C \left(|w_t''||J_t-J|+|w_t'||\nabla(J_t-J)|+|w_t|\nabla^2|J_t-J|\right)\\
&\leq& C\varepsilon_t.
\end{array}
\end{displaymath}

If $\dim M\geq 5$, then, since $M$ is compact, $(M,H,J)$ is embeddable. Let us consider an ACH manifold $(X,g)$ with CR infinity $(M,H,J)$. Let $J_X$ be a complex structure on $\overline{X}$ and let $z=(z_1,\ldots,z_{n+1})$ be complex coordinates near $\lbrace p_1,p_2\rbrace$. Then, by the normal form theorem of Chern and Moser, there is a boundary defining function $r$ on $\overline{X}$ such that $$r(z)=r_0(z)+\underset{1\leq j\leq n}{O}(|z_j|^4),$$ where $r_0(z)=\re(z_{n+1})-\frac{1}{4}\sum_{1\leq j\leq n}|z_j|^2$ is a boundary defining function for the %spherical half-space model in the Siegel domain
Heisenberg group \cite{CM74}. We glue the defining functions as follows: $$r_t=(1-\chi_t) r_0+\chi_t r.$$ The corresponding contact form is given by $\theta_t=i\left(\overline{\partial}-\partial\right)r_t$. The induced complex structure $J_t$ on $M$ is then given by the relation $d\theta(\cdot,J_t\cdot)=d\theta_t(\cdot,i\cdot)$. By construction, $J_t$ is spherical inside $U_t$ and coincides with $J$ outside $U'_t$, and $(J_t)$ $C^0$-converges to $J$. Moreover, since $|r_t-r|=O(x^4)$ and $|\nabla(r_t-r)|=O(x^3)$, we have, for some constant $C$, %$\mathrm{Scal}_W(J_t,\theta)\underset{t\rightarrow 0}{\longrightarrow}\mathrm{Scal}_W(J_0,\theta)$. 
% where $\tau=\iota_T T^\nabla$ is the pseudohermitian torsion, 
$$|\mathrm{Scal}_W(J_t,\theta)-\mathrm{Scal}_W(J,\theta)|\leq C \left(|w_t'|\left(|\nabla^2(r_t-r)|+|\nabla^3(r_t-r)|\right)+|w_t''||\nabla^2(r_t-r)|\right)\leq C\varepsilon_t.$$

%Moreover, 
%Let $U'\subset U$ be a smaller neighbourhood and $\hat{J}$ be a strictly pseudoconvex complex structure on $(M,H,\theta)$ coinciding with $J$ outside $U$ and such that $\hat{J}=\Phi^* J^c_\lambda$ on $U'$.
%by Lemma 3.4.3. in \cite{BR09}, there exist constants $C_p$ independent of $\varepsilon_t$ such that, on $U'_t$, $$\|J_t-J\|\leq C_0 \varepsilon_t^\frac{1}{2}\quad\text{and}\quad\|\nabla^p J_t-\nabla^p J\|\leq C_p \varepsilon_t^\frac{p}{2}\quad\text{for}\ p\geq 1,$$ where the norm is taken with respect to the standard metric induced by the embedding in $\C^n$. The family $(J_t)$ is thus $C^2$-convergent.
%We use a construction of a complex structure on an SPCR manifold $(M,H,J)$ which is spherical near the surgery locus $S$ and which coincides with $J$ elsewhere, due to Biquard and Rollin \cite{BR09}. Let $(z_i)$ denote complex coordinates near the surgery locus and $\varphi$ denote a K\"{a}hler potential compatible with $(M,H,J)$. We can suppose, using the normal form of Chern and Moser \cite{CM74}, that $$\varphi = \varphi_\Heis+O(|z|^4).$$ Let $I=(r_0,r_1)\subset\R_+^*$, let $u_I$ be a cut-off function
%\begin{displaymath}
%u_I(x)=\left\lbrace
%\begin{array}{ll}
%0 \ \mathrm{if} \ x<r_0,\\
%1 \ \mathrm{if} \ x>r_1,\\
%\end{array} \right.
%\end{displaymath}
%and $$\chi_I (z)=u_I\left(\left|\frac{z_n}{z_0}\right|\right).$$ Then the gluing $$\varphi_I = (1-\chi_I)\varphi_\Heis+\chi_I\varphi$$ is the K\"{a}hler potential of a manifold $(M,\tilde{H},\tilde{J})$ spherical near $S$ and coinciding with $(M,H,J)$ elsewhere.
\end{proof}

\begin{ex} If $(M,H)=(\Sph^{2n+1},H_0)$, using Theorem \ref{jl} we thus have the equality $$\sigma_c(\Sph^{1}\times\Sph^{2n},\tilde{H}_0)=\sigma_c(\Sph^{2n+1},H_0).$$
\end{ex}

\subsection{Disjoint union}

%Following \cite{Kob87}, Theorem \ref{CRkob} is obtained by the following inequality, proved in section \ref{glu}:

In the case of a connected sum, Theorem \ref{CRkobh} can be written the following way:

\begin{theo}\label{condis}
Let $(M_1,H_1)$ and $(M_2,H_2)$ be two compact SPC manifolds of dimension $2n+1$. Let $(M_1,H_1)\# (M_2,H_2)$ be their SPC connected sum, then $$\sigma_c\left((M_1,H_1)\# (M_2,H_2)\right)\geq \sigma_c\left((M_1,H_1)\sqcup (M_2,H_2)\right).$$
\end{theo}

\noindent Alongside with the hereunder computation of the right-hand side, this gives Theorem \ref{CRkob}.%$\sigma_c$ for a disjoint union of SPC manifolds.

%The disjoint union case is directly inspired from \cite{Kob87}.

\begin{prop}
\label{eqlem}
Let $(M_1,H_1)$ and $(M_2,H_2)$ be two compact SPC manifolds of dimension $2n+1$. Then
\begin{displaymath}
\sigma_c\left((M_1,H_1)\sqcup (M_2,H_2)\right)= \left\lbrace
\begin{array}{ll}
 -\left(|\sigma_c(M_1,H_1)|^{n+1}+|\sigma_c(M_2,H_2)|^{n+1}\right)^\frac{1}{n+1} & \mathrm{if} \ \sigma_c(M_1,H_1)\leq 0\\
 & \mathrm{and} \ \sigma_c(M_2,H_2)\leq 0,\\
 \min\left(\sigma_c(M_1,H_1),\sigma_c(M_2,H_2)\right) & \mathrm{otherwise}.
\end{array}\right.
\end{displaymath}
\end{prop}

\begin{proof}
Let us consider a unit volume strictly convex pseudohermitian structure $(J,\theta)$ on $(M_1,H_1)\sqcup (M_2,H_2)$. Let us denote, for $i$ in $\lbrace 1,2\rbrace,\ J_i:=J\mid_{M_i}$, $\theta_i:=\theta\mid_{M_i}$, $Y_i=Y_{CR}(M_i,H_i,J_i)$, and $\lambda_i$ in $\R_+^*$ which verifies $$\mathrm{Vol}(M_i,\lambda_i\theta_i)=\lambda_i^{n+1}\mathrm{Vol}(M_i,\theta_i)=1.$$ We recall that $S_W$ denotes the integral Webster scalar curvature. Since, for $i$ in $\lbrace 1,2\rbrace$, $\mathrm{Scal}_W(J,\lambda_i\theta_i)=\lambda_i^{-1}\mathrm{Scal}_W(J,\theta_i)$, we have
\begin{displaymath}
\begin{array}{cl}
\displaystyle{S_W(M_1\sqcup M_2,J,\theta)} & = \displaystyle{S_W(M_1,J_1,\theta_1)+S_W(M_2,J_2,\theta_2)}\\
& \displaystyle{= \lambda_1^{-n}S_W(M_1,J_1,\lambda_1\theta_1)+\lambda_2^{-n}S_W(M_2,J_2,\lambda_2\theta_2)}\\
& \displaystyle{= \mathrm{Vol}(M_1,\theta_1)^\frac{n}{n+1}S_W(M_1,J_1,\lambda_1\theta_1)+\mathrm{Vol}(M_2,\theta_2)^\frac{n}{n+1}S_W(M_2,J_2,\lambda_2\theta_2)}\\
& \displaystyle{\geq \mathrm{Vol}(M_1,\theta_1)^\frac{n}{n+1}Y_1+\mathrm{Vol}(M_2,\theta_2)^\frac{n}{n+1}Y_2},
\end{array}
\end{displaymath}
with equality when $\lambda_1\theta_1$ and $\lambda_2\theta_2$ are Yamabe contact forms on $(M_1,H_1,J_1)$ and $(M_2,H_2,J_2)$ respectively. Optimizing the right-hand side under the constraint $\mathrm{Vol}(M_1,\theta_1)+\mathrm{Vol}(M_2,\theta_2)=1$ yields
\begin{displaymath}
S_W(M_1\sqcup M_2,J,\theta)\geq \left\lbrace
\begin{array}{ll}
 -\left(|Y_1|^{n+1}+|Y_2|^{n+1}\right)^\frac{1}{n+1} & \mathrm{if}  \ Y_1\leq 0\ \mathrm{and} \ Y_2\leq 0,\\
 \min\left(Y_1,Y_2\right) & \mathrm{otherwise},
\end{array}\right.
\end{displaymath}
with equality when $\lambda_1\theta_1$ and $\lambda_2\theta_2$ are Yamabe contact forms and, in the first case, when $\frac{1}{\mathrm{Vol}(M_1,\theta_1)}=1+\left(\frac{Y_2}{Y_1}\right)^{n+1}$, and, in the second case, at the limit $\mathrm{Vol}(M_i,\theta_i)\rightarrow 0$, where $i\in\lbrace 1,2\rbrace$ verifies $Y_i=\max\left(Y_1,Y_2\right).$ Consequently,
%$$Y_{CR}\left((M_1,H_1,J_1)\sqcup (M_2,H_2,J_2)\right)=-\left(|Y_{CR}(M_1,H_1,J_1)|^{n+1}+|Y_{CR}(M_2,H_2,J_2)|^{n+1}\right)^\frac{1}{n+1}$$ if $Y_{CR}(M_1,H_1,J_1)\leq 0$ and $Y_{CR}(M_2,H_2,J_2)\leq 0$, $\min\left(Y_{CR}(M_1,H_1,J_1),Y_{CR}(M_2,H_2,J_2)\right)$ otherwise.
\begin{displaymath}
Y_{CR}\left((M_1,H_1,J_1)\sqcup (M_2,H_2,J_2)\right)= \left\lbrace
\begin{array}{ll}
 -\left(|Y_1|^{n+1}+|Y_2|^{n+1}\right)^\frac{1}{n+1} & \mathrm{if}  \ Y_1\leq 0\ \mathrm{and} \ Y_2\leq 0,\\
 \min\left(Y_1,Y_2\right) & \mathrm{otherwise},
\end{array}\right.
\end{displaymath}
hence the result.
\end{proof}

\section{A CR Gauss-Bonnet-LeBrun formula}
\label{comp}

We prove in this part Theorem \ref{CRleb}. Let us first recall some facts on the Burns-Epstein invariant of a CR manifold. Let $(M,H,J)$ be a compact SPCR $3$-manifold. The \emph{Burns-Epstein invariant} $\mu(M,H,J)$ is defined as the evaluation of a well-chosen de Rham cohomology class on the fundamental class $[M]$ in $H_3(M,\R)$ \cite{BE88}. In particular, we have the following estimates.%, which will be useful in Chapter \ref{Chap:4}.

\begin{prop}[\cite{Mar15}]
\label{BE}
The Burns-Epstein invariant of a compact SPCR $3$-manifold $(M,H,J)$ admitting a pseudo-Einstein contact form $\theta$ is given by $$\mu(M,H,J) = \frac{1}{4\pi^2}\int_M\left(|\tau(J,\theta)|^2-\frac{1}{4}\mathrm{Scal}_W(J,\theta)^2\right)\theta\wedge d\theta.$$% where $\tau$ is the pseudohermitian torsion.
\end{prop}

\begin{prop}[\cite{BE88}]
\label{BE2}
The Burns-Epstein invariant value of a circle bundle $(M,H,J)$ over a Riemann surface $\Sigma$ is $$\mu(M,H,J) = -\frac{|\chi(\Sigma)|}{4}+\frac{1}{12\pi}\int_\Sigma\frac{\left(\Delta\log |\mathrm{Scal}_W(J,\theta_J)|\right)^2}{|\mathrm{Scal}_W(J,\theta_J)|} d\mathrm{vol}_\Sigma,$$ where $\theta_J$ is the unique normal contact form on $(M,H,J)$ and $\chi(\Sigma)$ is the Euler characteristic of $\Sigma$.
\end{prop}

We now prove the CR analogue of a result due to C. LeBrun \cite{LeB99}.

%Theorem \ref{CRleb} is a consequence of Theorem \ref{CRleb2} below. Indeed, compact normal CR 3-manifolds are either quotients of $\Sph^3$, which have positive $\sigma_c$, or of a circle bundle $(X,H,J)$ over a Riemann surface $\Sigma$ of positive genus \cite{Bel01}. The Burns-Epstein invariant value for such a bundle is $$\mu(X,H,J) = -\frac{|\chi(\Sigma)|}{4}+\frac{1}{12\pi}\int_\Sigma\frac{\left(\Delta\log |\mathrm{Scal}_W(J,\theta_J)|\right)^2}{|\mathrm{Scal}_W(J,\theta_J)|} d\mathrm{vol}_\Sigma,$$ where $\theta_J$ is the unique normal contact form on $(X,H,J)$ and $\chi(\Sigma)$ is the Euler characteristic of $\Sigma$ \cite{BE88}. If $(X,H,J,\theta_{J})$ is also pseudo-Einstein, it has constant Webster scalar curvature. In this case, we have consequently $$\mu(X,H,J)=-\frac{|\chi(\Sigma)|}{4}.$$% hence Theorem \ref{CRleb}.

%\begin{theo}
%\label{CRleb2}
%If a compact SPC $3$-manifold $(M,H)$ admits for all $J\in\mathcal{J}$ a normal, pseudo-Einstein, Yamabe contact form $\theta_J$, if $\sigma_c(M,H)\leq 0$, %if $\mathcal{J}\ni J\mapsto \mu(M,H,J)$ realizes its supremum at $J_0$, 
%%and if for all $J\in\mathcal{J}$ there exists a Yamabe contact form, 
%then $$\sigma_c(M,H)=-4\pi\sqrt{-\underset{J\in\mathcal{J}}{\sup}\mu(M,H,J)}.$$ %where $\mu(M,H,J_0)$ is the Burns-Epstein invariant of $(M,H,J_0)$.
%\end{theo}

%Theorem \ref{CRleb} is a direct consequence of the following proposition:

%We first prove the following lemma:

\begin{prop}
\label{leb}
Let $(M,H,J)$ be a compact SPCR manifold of dimension $2n+1$ admitting a Yamabe contact form. Then $$|Y_{CR}(M,H,J)|^{n+1} = \inf_{\hat{\theta}\in[\theta]}\int_M|\mathrm{Scal}_W(J,\hat{\theta})|^{n+1} \hat{\theta}\wedge d\hat{\theta}^n,$$ and the infimum is realized by Yamabe contact forms.% Moreover, when $\sigma_c(M,H)\leq 0$, $$|\sigma_c(M,H)|^{n+1}=\underset{J\in \mathcal{J},\hat{\theta}\in[\theta]}{\inf}\int_M|\mathrm{Scal}_W(J,\hat{\theta})|^{n+1} \hat{\theta}\wedge d\hat{\theta}^n.$$
\end{prop}

\begin{proof}
%We follow \cite{LeB99}.
By H\"{o}lder's inequality, for all $\hat{\theta}\in[\theta]$, $$\left(\int_M|\mathrm{Scal}_W(J,\hat{\theta})|^{n+1}\hat{\theta}\wedge d\hat{\theta}^n\right)^\frac{1}{n+1}\geq\frac{\int_M \mathrm{Scal}_W(J,\hat{\theta})\hat{\theta}\wedge d\hat{\theta}^n}{\left(\int_M\hat{\theta}\wedge d\hat{\theta}^n\right)^\frac{n}{n+1}},$$ with equality if and only if $\mathrm{Scal}_W(J,\hat{\theta})$ is a non-negative constant. If $Y_{CR}(M,H,J)\geq 0$, the claim follows from the fact that there exists a Yamabe contact form.

If $Y_{CR}(M,H,J)< 0$, let $\tilde{\theta}\in[\theta]$ be a Yamabe contact form. Let us consider $\hat{\theta}\in[\theta]$ and $u\in C^\infty(M, \R_+^*)$ such that $\hat{\theta}=u^\frac{2}{n}\tilde{\theta}$. Then $$\mathrm{Scal}_W(J,\hat{\theta})u^\frac{n+2}{n}=2\left(\frac{n+1}{n}\right)\Delta_b u+\mathrm{Scal}_W(J,\tilde{\theta})\cdot u,$$ so that
\begin{displaymath}
\begin{array}{ll}
\int_M \mathrm{Scal}_W(J,\hat{\theta})u^\frac{2}{n} \tilde{\theta}\wedge d\tilde{\theta}^n & = \int_M \left(2\left(\frac{n+1}{n}\right)\frac{\Delta_b u}{u}+\mathrm{Scal}_W(J,\tilde{\theta})\right) \tilde{\theta}\wedge d\tilde{\theta}^n\\
& = \int_M \left(-2\left(\frac{n+1}{n}\right)\frac{|d_bu|^2}{u^2}+\mathrm{Scal}_W(J,\tilde{\theta})\right) \tilde{\theta}\wedge d\tilde{\theta}^n\\
& \leq S_W(M,J,\tilde{\theta}),
\end{array}
\end{displaymath} and by H\"{o}lder's inequality,
\begin{displaymath}
\begin{array}{ll}
\left(\int_M |\mathrm{Scal}_W(J,\hat{\theta})|^{n+1} \hat{\theta}\wedge d\hat{\theta}^n\right)^\frac{1}{n+1}
 & = \left(\int_M |\mathrm{Scal}_W(J,\hat{\theta})u^\frac{2}{n}|^{n+1} \tilde{\theta}\wedge d\tilde{\theta}^n\right)^\frac{1}{n+1}\\
 & \geq -\frac{\int_M \mathrm{Scal}_W(J,\hat{\theta})u^\frac{2}{n}\tilde{\theta}\wedge d\tilde{\theta}^n}{\left(\int_M \tilde{\theta}\wedge d\tilde{\theta}^n\right)^\frac{n}{n+1}}\\
& \geq -\frac{S_W(M,J,\tilde{\theta})}{\left(\int_M \tilde{\theta}\wedge d\tilde{\theta}^n\right)^\frac{n}{n+1}}\\
& = |Y_{CR}(M,H,J)|,
\end{array}
\end{displaymath} with equality if and only if $u$ is a constant, which proves the desired equality.% From it, we have that $$\inf_{J\in\mathcal{J}}|Y_{CR}(M,H,J)|^{n+1} = \underset{J\in \mathcal{J},\hat{\theta}\in[\theta]}{\inf}\int_M|\mathrm{Scal}_W(J,\hat{\theta})|^{n+1} \hat{\theta}\wedge d\hat{\theta}^n.$$ If $\sigma_c(M,H)\leq 0$, the left-hand side of this equality is equal to $|\sigma_c(M,H)|^{n+1}$, hence the second equality.
\end{proof}

This proposition yields the following estimate on $Y_{CR}$, which implies Theorem \ref{CRleb}.
%
%\begin{cor}
%Let $(M,H,J)$ be a compact SPCR $3$-manifold admitting an Einstein, Yamabe contact form. If its Burns-Epstein invariant $\mu(M,H,J)$ is nonpositive, then $$|Y_{CR}(M,H,J)|= 4\pi\sqrt{-\mu(M,H,J)}.$$
%%Let $(M,H)$ be an SPC compact $3$-manifold admitting a pseudo-Einstein pseudohermitian structure $(J,\theta)$ such that $Y_{CR}(M,H,J)\geq 0$ and $\mu(M,H,J)\leq 0$. Then $$\sigma_c(M,H)\geq 4\pi\underset{J\in\mathcal{J}}{\sup}\sqrt{-\mu(M,H,J)}.$$
%\end{cor}
%
%\begin{proof}
%Let $\theta$ be an Einstein contact form on $(M,H,J)$. Then, by Proposition \ref{BE}, $$\int_M \mathrm{Scal}_W(J,\theta)^2\theta\wedge d\theta= -16\pi^2\mu(M,H,J).$$ If $\theta$ is also Yamabe, then by Proposition \ref{leb}, $$Y_{CR}(M,H,J)^2= -16\pi^2\mu(M,H,J).$$
%\end{proof}

\begin{cor}
Let $(M,H,J)$ be a circle bundle over a Riemann surface $\Sigma$ of positive genus admitting an Einstein contact form. Then $$Y_{CR}(M,H,J)= -2\pi\sqrt{-\chi(\Sigma)}.$$
%Let $(M,H)$ be an SPC compact $3$-manifold admitting a pseudo-Einstein pseudohermitian structure $(J,\theta)$ such that $Y_{CR}(M,H,J)\geq 0$ and $\mu(M,H,J)\leq 0$. Then $$\sigma_c(M,H)\geq 4\pi\underset{J\in\mathcal{J}}{\sup}\sqrt{-\mu(M,H,J)}.$$
\end{cor}

\begin{proof}
Let $\theta$ be an Einstein contact form on $(M,H,J)$. By Propositions \ref{BE} and \ref{BE2}, $$\int_M \mathrm{Scal}_W(J,\theta)^2\theta\wedge d\theta= -16\pi^2\mu(M,H,J)=4\pi^2|\chi(\Sigma)|.$$ Then by Proposition \ref{leb}, 
\begin{equation}\label{ineqyam}
Y_{CR}(M,H,J)^2\leq  4\pi^2|\chi(\Sigma)|.
\end{equation}
If $\Sigma$ is a torus, this implies that $Y_{CR}(M,H,J)=0$. Otherwise, $(M,H,J)$ admits a contact form of negative Webster scalar curvature, hence $Y_{CR}(M,H,J)\leq 0$ by Proposition \ref{pos}. In all cases, $Y_{CR}(M,H,J)\leq 0<Y_{CR}(\Sph^3,H_0,J_0).$ By Theorems \ref{jlb} and \ref{einwang}, $\theta$ is thus a Yamabe contact form, hence the inequality (\ref{ineqyam}) is an equality.
\end{proof}

\bibliographystyle{alpha}

\bibliography{../Biblio}

\end{document}